\documentclass[10pt]{article}
\usepackage{amsmath,amsthm,amsfonts,amssymb,amscd,latexsym,eucal}
\usepackage[a4paper,top=3cm,bottom=4.5cm,left=2.5cm,right=2.5cm]{geometry}
\usepackage{amsmath}
\usepackage{amsfonts}
\usepackage{amssymb}
\usepackage{graphicx}
\newtheorem{theorem}{Theorem}

\newtheorem{corollary}{Corollary}
\newtheorem{definition}{Definition}
\newtheorem{lemma}{Lemma}
\newtheorem{proposition}{Proposition}
\newtheorem{remark}{Remark}

\begin {document}

\author{Carlo Pandiscia}
\title{Covariant GNS Representation for C*-Dynamical Systems}
\date{}
\maketitle

\begin {abstract}
We extend the covariant GNS representation of Niculescu, Str\"{o}h  and Zsid\'{o} for C*-dynamical systems with time-evolution of the system (dynamics) a homomorphism of C*-algebras, to any dynamical systems, where the dynamics is an unital completely positive map. We give also an overview on its application to the reversible dilation theory as formulated by B. Kummerer. 
\end {abstract}

\section{Introduction}
Let $\mathcal{H}$ be a Hilbert space, with $\mathcal{B}(\mathcal{H})$ we denote the C*-algebra of all bounded linear operators on $\mathcal{H}$. Furthermore $wo$ stands for weak operator topology on $\mathcal{B}(\mathcal{H})$, $so$ for the strong operator topology, and $\omega$ for the weak topology defined on $\mathcal{B}( \mathcal{H})$ (see \cite{Stra-Zsi}). Furthermore a linear form $\omega $ on $\mathcal{B}( \mathcal{H})$ is said be \textit{normal} if is a $w$-continuous linear form and a homomorphism $\Phi :\mathfrak{A}\rightarrow \mathfrak{B}$ between unital C*-algebras is an unital *-multiplicative map, while a representation of a C*-algebra $\mathfrak{A}$ on the Hilbert space $\mathcal{H}$ is a homomorphism $\pi:\mathfrak{A}\rightarrow \mathcal{B}( \mathcal{H})$.
\newline
A C*-dynamical system is a triple $(\mathfrak{A},\Phi,\varphi)$ constitued by a C*-algebra with unit $\mathfrak{A}$, an unital completely positive map (briefly \textit{ucp-map}) $\Phi:\mathfrak{A}\rightarrow \mathfrak{A}$ and a state $\varphi$ on $\mathfrak{A}$ such that $\varphi\circ\Phi=\varphi$. Furthermore the ucp-map $\Phi$ is said be the dynamics of our C*-dynamical system.
\newline
In particular a C*-dynamical system $(\mathfrak{M},\Phi,\varphi)$ constituted by a von Neumann algebra $\mathfrak{M}$, normal ucp-map $\Phi$ and by a normal faithful state $\varphi$, will be called a W*-dynamical system.  

\smallskip
To enter the topic of this paper, let $(\mathcal{H}_{\varphi},\pi_{\varphi},\Omega_{\varphi})$ be the GNS representation of
$\varphi$, it is well know that there is an unique linear contraction
$\mathbf{U}_{\Phi,\varphi}$ of $\mathfrak{B}\left(\mathcal{H}_{\varphi
}\right)  $ such that, for any $a\in\mathfrak{A}$, we have
\begin{equation}
\mathbf{U}_{\Phi,\varphi}\pi_{\varphi}(a)\Omega_{\varphi}=
\pi_{\varphi}(\Phi(a))\Omega_{\varphi} \label{contrazione 1}.
\end{equation}
Moreover, it is simple to prove that if $\Phi$ is a homomorphism, then the
contraction $\mathbf{U}_{\Phi,\varphi}$ is an isometry on $\mathcal{H}
_{\varphi}$ and for any $a\in\mathfrak{A}$ we obtain
\begin{equation}
\mathbf{U}_{\Phi,\varphi}\pi_{\varphi}(a)=\pi_{\varphi}(\Phi(a))\mathbf{U}_{\Phi,\varphi}.
\end{equation}
If the support projection $s(\varphi)$ of $\varphi$ in the second dual $\mathfrak{A}^{\ast\ast}$ is central (this happens if and only if the vector $\Omega_{\varphi}$ is cyclic for $\pi_{\varphi}(\mathfrak{A})^{\prime}$) there exists a W*-dynamical system $(\pi_{\varphi}(\mathfrak{A})^{\prime\prime},\Phi_{\bullet},\varphi_{\bullet})$, where the dynamics
$\Phi_{\bullet}:\pi_{\varphi}(\mathfrak{A})^{\prime\prime}
\rightarrow\pi_{\varphi}(\mathfrak{A})^{\prime\prime}$ is the
normal ucp-map thus defined:
\begin{equation}
\Phi_{\bullet}(X)\Omega_{\varphi}=\mathbf{U}_{\Phi,\varphi}X\Omega_{\varphi} \ \ \text{for all} \ X\in\pi_{\varphi}(\mathfrak{A})^{\prime\prime},
\end{equation}
while $\varphi_{\bullet}$ is the normal faithful state
\begin{equation}
\varphi_{\bullet}(X)=\left\langle \Omega_{\varphi},
X\Omega_{\varphi}\right\rangle \ \ \text{for all} \ X\in\pi_{\varphi}(\mathfrak{A})^{\prime\prime}.
\end{equation}
In \cite{NSZ}, Niculescu, Str\"{o}h and Zsid\'{o}, using the minimal unitary dilation of the contraction $\mathbf{U}_{\Phi,\varphi}$ (see \cite{Nagy-Foias}), have proved the existence of a representation that generalizes the GNS representation associated to C*-dynamical system with dynamics $\Phi$ a homomorphism (i.e. C*-dynamical system with multiplicative dynamics), called \emph{the covariant GNS representation}, briefly CGNS representation. 
Specifically, they  proved the existence of a quadruple $(\mathcal{H},\pi ,\mathbf{U},\Omega) $ constituted by a unique, up to equivalence, representation $\pi :\mathfrak{A}\rightarrow \mathfrak{B}(\mathcal{H}) $, an unitary operator $\mathbf{U}$ on Hilbert space $\mathcal{H}$ and a vector $\Omega$ belonging to $\mathcal{H}$ such that
\begin{itemize}
\item[a)]
$\pi(\Phi(a))=\mathbf{U}\pi(a)\mathbf{U}^*$, for all $a\in\mathfrak{A}$;
\item[b)]
The pair $(\mathbf{U},\mathcal{H})$ is the minimal unitary dilation of $(\mathbf{U}_{\Phi,\varphi},\mathcal{H}_{\varphi})$ and $\mathbf{U}\Omega=\Omega$;
\item[c)]
$\Omega$ is a cyclic vector for the *-subalgebra $\mathfrak{B}$ generated by the set $\bigcup\limits_{k\in\mathbb{Z}}^{\infty}\{\mathbf{U}^{k}\pi(a)\mathbf{U}^{-k}:a\in\mathfrak{A}\}$;
\item[d)] 
$\varphi(a)=\left\langle {\Omega,\pi(a)\Omega}\right\rangle$ for all $a\in\mathfrak{A}$.
\end{itemize} 
\smallskip
In this paper we extend the CGNS representation, previously given for $C^*$-dynamical system with multiplicative dynamics $\Phi$, to a generic C*-dynamical system i.e. having only an ucp-map.
\newline
The existence of a such CGNS representation easily prove that the W*-dynamical system associated to a C*-dynamical system with multiplicative dynamics, admits a reversible dilation in the direction of Kummerer in \cite {Kummerer}.
\newline
After a summary of previous notation, the paper is organized as follow.
\newline
In section \ref{COVARIANT-GNS} using the Stinespring representation of completely positive maps and the inductive limit of directed systems of Hilbert space , we construct the CGNS representation of a dynamical system.
in section \ref{REVERSIBLE-DIL} we will show that the minimal reversible dilation of the W*-dynamical system associated to C*-dynamical systems with multiplicative dynamics, satisfies the ergodic properties of the original dynamical system. 
\newline
Furthermore we shall prove that a C*-dynamical system which admits a right inverse, i.e.\ an ucp-map $\Psi:\mathfrak{A}\rightarrow\mathfrak{A}$ such that for each $a$ belong to $\mathfrak{A}$ we have $\Phi(\Psi(a)))=a$, its associated W*-dynamical system $(\pi_{\varphi}(\mathfrak{A})^{\prime\prime},\Phi_{\bullet},\varphi_{\bullet})$ admit a minimal reversible dilation.

\section{Covariant GNS associated to the dynamical systems} \label{COVARIANT-GNS}
Before getting into a discussion of the covariant GNS representation, let us recall briefly the Stinespring's theorem of a completely positive maps (see \cite{Paulsen}).
\newline
We consider a Hilbert space $\mathcal{H}$ and a C*-subalgebra with unit  $\mathfrak{A}$ of  $\mathcal{B}(\mathcal{H})$, the Stinespring representation associated to an ucp-map $\Phi:\mathfrak{A}\rightarrow\mathfrak{A}$ is a triple $(\mathbf{V}_{\Phi},\sigma_{\Phi},\mathcal{L}_{\Phi})$, constituted by a Hilbert space $\mathcal{L}_{\Phi}$, a representation $\sigma_{\Phi}:\mathfrak{A}\rightarrow\mathcal{B}(\mathcal{L}_{\Phi}) $ 
and a linear contraction $\mathbf{V}_{\Phi}:\mathcal{H}\rightarrow\mathcal{L}_{\Phi}$ 
such that for $a\in\mathfrak{A}$ we have
\begin{equation}
\Phi(a)=\mathbf{V}_{\Phi}^*\sigma_{\Phi}(a)\mathbf{V}_{\Phi}.
\label{stine}
\end{equation}
We recall that on the algebraic tensor 
$\mathfrak{A}\otimes\mathcal{H}$ 
we can define a semi-inner product by
\[
\left\langle a_{1}\otimes\Psi_{1},a_{2}\otimes\Psi_{2}\right\rangle_{\Phi}=\left\langle \Psi_{1},\Phi\left( a_{1}^{\ast}a_{2}\right)\Psi_{2}\right\rangle_{\mathcal{H}},
\]
for all 
$a_{1},a_{2}\in\mathfrak{A}$ 
and 
$\Psi_{1},\Psi_{2}\in\mathcal{H}$
furthermore the Hilbert space 
$\mathcal{L}_{\Phi}$ 
is the completion of the quotient space 
$\mathfrak{A}\overline{\otimes}_{\Phi}\mathcal{H}$
 of 
$\mathfrak{A}\otimes\mathcal{H}$ 
by the linear subspace
\[
\left\{X\in\mathfrak{A}\otimes\mathcal{H}:\left\langle X,X\right\rangle
_{\Phi}=0\right\}
\]
with inner product induced by 
$\left\langle \cdot \ ,\cdot\right\rangle _{\Phi}$. 
We shall denote the image at 
$a\otimes\Psi\in\mathfrak{A}\otimes\mathcal{H}$ 
in 
$\mathfrak{A}\overline{\otimes}_{\Phi}\mathcal{H}$ 
by 
$a\overline{\otimes}_{\Phi}\Psi,$ 
so that we have
\[
\left\langle a_{1}\overline{\otimes}_{\Phi}\Psi_{2},a_{2}\overline{\otimes
}_{\Phi}\Psi_{2}\right\rangle _{\mathcal{L}_{\Phi}}=\left\langle \Psi_{1}
,\Phi\left(a_{1}^{\ast}a_{2}\right)\Psi_{2}\right\rangle _{\mathcal{H}},
\]
for all $a_{1},a_{2}\in\mathfrak{A}$ and $\Psi_{1},\Psi_{2}
\in\mathcal{H}$.
\newline
Moreover
$\sigma_{\Phi}\left(a\right)\left(x\overline{\otimes}_{\Phi}\Psi\right)  =ax\otimes_{\Phi}\Psi,$ 
for each 
$x\overline{\otimes}_{\Phi}\Psi\in\mathcal{L}_{\Phi}$ 
and 
$\mathbf{V}_{\Phi}\Psi=\mathbf{1}\overline{\otimes}_{\Phi}\Psi$ 
for each 
$\Psi\in\mathcal{H}.
$\newline 
Since $\Phi$ is unital map, the linear operator 
$\mathbf{V}_{\Phi}$ 
is an isometry with adjoint 
$\mathbf{V}_{\Phi}^{\ast}$ 
defined by
\[
\mathbf{V}_{\Phi}^*a\overline{\otimes}_{\Phi}\Psi=\Phi(a)\Psi,
\] 
for all $a\in\mathfrak{A}$ and $\Psi\in\mathcal{H}$.
\newline
We recall that the multiplicative domain of the ucp-map
$\Phi:\mathfrak{A}\rightarrow\mathfrak{A}$ is the set such defined:
\[
\mathcal{D}_{\Phi}=\{a\in\mathfrak{A}:\Phi(a^*)\Phi(a)=\Phi(a^*a) \  \text{and} \  \Phi(a)\Phi(a^*)=\Phi(aa^*) \}.
\]
We have the following implications (See \cite{Paulsen}):
\newline  
An element $ a\in\mathcal{D}_{\Phi}$   if and only if $\Phi(a)\Phi(x)=\Phi(ax)$  and   $\Phi(x)\Phi(a)=\Phi(xa)$ for all $x\in\mathfrak{A}$.
\newline
Then the set $\mathcal{D}_{\Phi}$ is an unital C*-subalgebra of $\mathfrak{A}$.
\begin{proposition}
For any $x\in\mathcal{D}_{\Phi}$ we have:
\[
\sigma_{\Phi}\left(x\right)\mathbf{V}_{\Phi}\mathbf{V}_{\Phi}^{\ast}=
\mathbf{V}_{\Phi}\mathbf{V}_{\Phi}^{\ast}\sigma_{\Phi}\left(x\right),
\]
it follows that $\Phi$ is homomorphism if and only if $\mathbf{V}_{\Phi}$ is an unitary.

\end{proposition}

\begin {proof}
For each $\Psi\in\mathcal{H}$ we obtain the following implications:
\[
a\overline{\otimes}_{\Phi}\Psi=\mathbf{1}\overline{\otimes}_{\Phi}\Phi\left(
a\right)\Psi\ \ \ \text{if and only if} \ \ \ \Phi\left(a^{\ast}a\right)
=\Phi\left( a^{\ast}\right)\Phi\left(a\right)  ,
\]
since
\[
\left\Vert a\overline{\otimes}_{\Phi}\Psi-1\overline{\otimes}_{\Phi}
\Phi\left(a\right)\Psi\right\Vert =\left\langle \Psi,\Phi\left(  a^{\ast
}a\right)\Psi\right\rangle -\left\langle \Psi,\Phi\left( a^{\ast}\right)
\Phi\left(a\right)\Psi\right\rangle .
\]
Furthermore, for each 
$a\in\mathfrak{A}$ 
and 
$\Psi\in\mathcal{H}$ 
we have
$\mathbf{V}_{\Phi}\mathbf{V}_{\Phi}^{\ast}a\overline{\otimes}_{\Phi}
\Psi=\mathbf{1}\overline{\otimes}_{\Phi}\Phi\left(a\right)\Psi.$
\end {proof}
Let $(\mathfrak{A},\Phi,\varphi)$ be an any C*-dynamical system, we set with $(\mathcal{L}_{1},\sigma_{1},\mathbf{V}_{0})  $ the
Stinespring representation of the normal
ucp-map $\Phi_{0}:\mathfrak{A\rightarrow B}(\mathcal{H}_{\varphi})$ defined by
\[
\Phi_{0}(a)=\pi_{\varphi}(\Phi(a)) \ \ \text{for all} \ a\in\mathfrak{A}. 
\]
The $\sigma_{1}:\mathfrak{A}\rightarrow\mathfrak{B}(\mathcal{L}_{1})$ is a representation on the Hilbert space 
$\mathcal{L}_{1}=\mathfrak{A}\overline{\otimes}_{\Phi_{0}}\mathcal{H}_{\varphi}$ 
such that:
\[
\Phi_{0}\left(  a\right)  =\mathbf{V}_{0}^{\ast}\sigma_{1}(a)\mathbf{V}_{0} \ \ \text{for all} \  a\in\mathfrak{A}, 
\]
with $\mathbf{V}_{0}:\mathcal{H}_{\varphi}\rightarrow\mathcal{L}_{1}$ linear
isometry thus defined
\[
\mathbf{V}_{0}h=1\overline{\otimes}_{\Phi_{0}}h \ \ \text{for all} \ h\in\mathcal{H}_{\varphi}.
\]
We recall that if $a_{1}\overline{\otimes}_{\Phi_{0}}h_{1}$ and $a_{2}
\overline{\otimes}_{\Phi_{0}}h_{2}$ are elements of $\mathcal{L}_{1}$ we have:
\[
\left\langle a_{1}\overline{\otimes}_{\Phi_{0}}h_{1},a_{2}\overline{\otimes}_{\Phi}h_{2}\right\rangle _{\mathcal{L}_{1}}=\left\langle h_{1},\Phi
_{0}\left(  a_{1}^{\ast}a_{2}\right)  h_{2}\right\rangle _{\mathcal{H}
_{\varphi}},
\]
furthermore for each $a\overline{\otimes}_{\Phi_{0}}h\in\mathcal{L}_{1}$
\[
\mathbf{V}_{0}^{\ast}a\overline{\otimes}_{\Phi_{0}}h=\Phi_{0}(a)
h.
\]
We have the follow lemma:

\begin{lemma}
\label{lemma-fattor}
There exists a linear isometry $\Lambda_{0}:\mathcal{H}_{\varphi}\rightarrow\mathcal{L}_{1}$ such that for any $a\in\mathfrak{A}$ we have
\begin{equation}
\Lambda_{0}\pi_{\varphi}(a)\Omega_{\varphi}=a\overline{\otimes
}_{\Phi_{0}}\Omega_{\varphi},
\label{lambda 0}%
\end{equation}
and
\[
\sigma_{1}\left(a\right)  \Lambda_{0}=\Lambda_{0}\pi_{\varphi}\left(a\right).
\]
Moreover the linear contraction $\mathbf{U}_{\Phi,\varphi}$ of $\mathfrak{B}
(\mathcal{H}_{\varphi})$ has the following factorization:
\begin{equation}
\mathbf{U}_{\Phi,\varphi}=\mathbf{V}_{0}^{\ast}\Lambda_{0}.
\label{fattorizzazione 0}%
\end{equation}
\end{lemma}

\begin {proof}
For any $a\in\mathfrak{A}$ we have
\[
\left\Vert a\overline{\otimes}_{\Phi}\Omega_{\varphi}\right\Vert
^{2}=\left\langle a\overline{\otimes}_{\Phi}\Omega_{\varphi},a\overline
{\otimes}_{\Phi}\Omega_{\varphi}\right\rangle _{\mathcal{L}_{1}}=\left\langle
\Omega_{\varphi},\Phi_{0}\left(a^{\ast}a\right)\Omega_{\varphi
}\right\rangle _{\mathcal{H}_{\varphi}}=\varphi\left(a^{\ast}a\right)
=\left\Vert \pi_{\varphi}\left(a\right)\Omega_{\varphi}\right\Vert^{2}.
\]
Then the linear map $\Lambda_{0}:\mathcal{H}_{\varphi}\rightarrow\mathcal{L}_{1}$ defined by the relationship \ref{lambda 0} it is well defined and isometric, follows that can be extended to all Hilbert space
$\mathcal{H}_{\varphi}$. Furthermore for each $x\in\mathfrak{A}$ we obtain:
\[
\sigma_{1}\left(a\right)\Lambda_{0}\pi_{\varphi}\left(x\right)
\Omega_{\varphi}=\sigma_{1}\left(a\right)x\overline{\otimes}_{\Phi}
\Omega_{\varphi}=ax\overline{\otimes}_{\Phi}\Omega_{\varphi}=\Lambda_{0}
\pi_{\varphi}\left(  ax\right)\Omega_{\varphi}=\Lambda_{0}\pi_{\varphi
}\left(  a\right)\pi_{\varphi}\left(x\right)\Omega_{\varphi},
\]
and
\[
\mathbf{V}_{0}^{\ast}\Lambda_{0}\pi_{\varphi}\left(  x\right)  \Omega
_{\varphi}=\mathbf{V}_{0}^{\ast}x\overline{\otimes}_{\Phi}\Omega_{\varphi
}=\Phi_{0}\left(  x\right)  \Omega_{\varphi}=\mathbf{U}_{\Phi,\varphi}%
\pi_{\varphi}\left(  x\right)  \Omega_{\varphi}.
\]
\end {proof}
We consider the normal ucp-map $\Phi_{1}:\mathfrak{A}\rightarrow
\mathfrak{B}(\mathcal{L}_{1})$ defined by
\[
\Phi_{1}(a)=\sigma_{1}(\Phi(a)) \ \ \text{for all} \ a\in\mathfrak{A} 
\]
and its Stinespring representation $(\mathcal{L}_{2},\sigma
_{2},\mathbf{V}_{1})$, with $\mathcal{L}_{2}=\mathfrak{A}
\overline{\otimes}_{\Phi_{1}}\mathcal{L}_{1}$ and $\mathbf{V}_{1}
:\mathcal{L}_{1}\rightarrow\mathcal{L}_{2}$ 
where for each $a\in\mathfrak{A}$ we have 
\[
\Phi_{1}(a)=\mathbf{V}_{1}^*\sigma_{2}(a)\mathbf{V}_{1}.
\]
We define a linear isometry $\Lambda_{1}:\mathcal{L}_{1}\rightarrow
\mathcal{L}_{2}$ as follows:
\[
\Lambda_{1}\sum\limits_{i}^{n}a_{i}\overline{\otimes}_{\Phi_{0}}h_{i}=
\sum \limits_{i}^{n}a_{i}\overline{\otimes}_{\Phi_{1}}\Lambda_{0}h_{i},
\]
for all $a_{i}\in\mathfrak{A}$ and $h_{i}\in\mathcal{H}_{\varphi}$, for each $i=1,2...n$.
\newline 
The $\Lambda_{1}$ is well defined operator, since for each $i,j$ we have:
\[
\left\langle a_{i}\overline{\otimes}_{\Phi_{1}}\Lambda_{0}h_{i},a_{j}
\overline{\otimes}_{\Phi_{1}}\Lambda_{0}h_{j}\right\rangle =\left\langle
\Lambda_{0}h_{i},\Phi_{1}(a_{i}^*a_{j})\Lambda_{0}h_{j}\right\rangle =\left\langle h_{i},\Lambda_{0}^*\Phi_{1}(a_{i}^*a_{j})\Lambda_{0}h_{j}\right\rangle=
\]
\[
=\left\langle h_{i},\Lambda_{0}^*\sigma_{1}(\Phi(a_{i}^*a_{j}))\Lambda_{0}h_{j}\right\rangle =
\left\langle h_{i},\pi_{\varphi}(\Phi(a_{i}^*a_{j}))\Lambda_{0}h_{j}\right\rangle =
\left\langle a_{i}\overline{\otimes}_{\Phi_{0}}h_{i},a_{j}\overline{\otimes}_{\Phi_{0}}h_{j}\right\rangle.
\] 
It is  simple to prove that for any $a\in\mathfrak{A}$, we have:
\[
\sigma_{2}(a)\Lambda_{1}=\Lambda_{1}\sigma_{1}(a)
\text{ \ \ \ and \ \ \ \ }\Lambda_{1}^*\sigma_{2}(a)\Lambda_{1}=\sigma_{1}(a).
\]
Furthermore the diagram
\[
\begin{array}
[c]{ccccc}%
\mathcal{H}_{\varphi} & \overset{\Lambda_{0}}{\longrightarrow} &
\mathcal{L}_{1} & \overset{\Lambda_{1}}{\longrightarrow} & \mathcal{L}_{2}\\
\downarrow & \ \ \overset{\mathbf{V}_{0}}{\searrow} & \downarrow &
\ \ \overset{\mathbf{V}_{1}}{\searrow} & \downarrow\\
\mathcal{H}_{\varphi} & \overset{\Lambda_{0}}{\longrightarrow} &
\mathcal{L}_{1} & \overset{\Lambda_{1}}{\longrightarrow} & \mathcal{L}_{2}%
\end{array}
\]
is commutative
\[
\mathbf{V}_{1}\Lambda_{0}=\Lambda_{1}\mathbf{V}_{0},
\]
with
\[
\Lambda_{0}\mathbf{V}_{0}^{\ast}=\mathbf{V}_{1}^{\ast}\Lambda_{1}.
\]
In fact for each $a\overline{\otimes}_{\Phi_{0}}h$ belong to $\mathcal{L}_{1}$ we have:
\[
\mathbf{V}_{1}^*\Lambda_{1}a\overline{\otimes}_{\Phi_{0}}h=\mathbf{V}_{1}^* a\overline{\otimes}_{\Phi_{1}}\Lambda_{0}h=\Phi_{1}(a)\Lambda_{0}h=\sigma_{1}(\Phi(a))\Lambda_{0}h=
\Lambda_{0}\pi_{\varphi}(\Phi(a))h=\Lambda_{0}\mathbf{V}_{0}^* a\overline{\otimes}_{\Phi_{0}}h.
\]
Iterating the procedure, for every natural number $n$, we have the normal ucp-map $\Phi_{n}:\mathfrak{A}
\rightarrow\mathfrak{B}(\mathcal{L}_{n})$ defined by:
\begin{equation}
\Phi_{n}(a)=\sigma_{n}(\Phi(a)) \ \ \text{for all} \ a\in\mathfrak{A}.
\label{ucpn-map} 
\end{equation}
We set again with $(\mathcal{L}_{n+1},\sigma_{n+1},\mathbf{V}_{n})$ its Stinespring representation.
Then we have the Hilbert space $\mathcal{L}_{n+1}=\mathfrak{A}\overline{\otimes}_{\Phi_{n}}\mathcal{L}_{n}$, the *-representation $\sigma_{n+1}:\mathfrak{A}\rightarrow\mathfrak{B}(\mathcal{L}_{n+1})$ and the linear isometry $\mathbf{V}_{n}:\mathcal{L}_{n}\rightarrow\mathcal{L}_{n+1}$, such that for each $a\in\mathfrak{A}$ we have
\[
\Phi_{n}(a)=\mathbf{V}_{n}^*\sigma_{n+1}(a)
,\mathbf{V}_{n}.
\]
Moreover we obtain a linear operator $\Lambda_{n}:\mathcal{L}_{n}\rightarrow
\mathcal{L}_{n+1}$ thus defined:
\[
\Lambda_{n}a\overline{\otimes}_{\Phi_{n-1}}\psi=a\overline{\otimes}_{\Phi_{n}%
}\Lambda_{n-1}\psi,
\]
for all $a\in\mathfrak{A},$ and $\psi\in\mathcal{L}_{n-1}$.
\newline 
The operator $\Lambda_{n}$ is an isometry furthermore for each natural number $n$ and
$a\in\mathfrak{A}$ we obtain
\begin{itemize}
\item[(a)] $\sigma_{n}\left(  a\right)  \Lambda_{n-1}=\Lambda_{n-1}\sigma
_{n-1}\left(  a\right)  ;$

\item[(b)] $\Lambda_{n}^{\ast}\sigma_{n}\left(  a\right)  \Lambda_{n}=\sigma
_{n-1}\left(  a\right)  ;$

\item[(c)] $\mathbf{V}_{n}\Lambda_{n-1}=\Lambda_{n}\mathbf{V}_{n-1};$

\item[(d)] $\Lambda_{n-1}\mathbf{V}_{n-1}^{\ast}=\mathbf{V}_{n}^{\ast}\Lambda_{n}$,
\end{itemize}
hence we have the commutative diagram
\begin{equation}
\begin{array}
[c]{ccccccccc}%
... & \longrightarrow & \mathcal{L}_{n-1} & \overset{\Lambda_{n-1}}{\longrightarrow} & \mathcal{L}_{n} & \overset{\Lambda_{n}}{\longrightarrow}
& \mathcal{L}_{n+1} & \overset{\Lambda_{n+1}}{\longrightarrow} & ....\\
& \sigma_{n-1}\left(  a\right)  & \downarrow & \overset{\mathbf{V}_{n-1}
}{\searrow} & \sigma_{n}\left(  a\right)  \downarrow & \overset{\mathbf{V}
_{n}}{\searrow} & \sigma_{n+1}(a)  \downarrow &
\overset{\mathbf{V}_{n+1}}{\searrow} & \\
... & \longrightarrow & \mathcal{L}_{n-1} & \overset{\Lambda_{n-1}
}{\longrightarrow} & \mathcal{L}_{n} & \overset{\Lambda_{n}}{\longrightarrow}
& \mathcal{L}_{n+1} & \overset{\Lambda_{n+1}}{\longrightarrow} & ....
\end{array}
\label{commutative diagram}
\end{equation}
We have a directed system of Hilbert spaces $(\mathcal{L}_{n},\Xi_{n,m})$ where the isometries $\Xi_{n,m}:\mathcal{L}_{m}\rightarrow\mathcal{L}_{n}$ 
for $m\leq n$, $\ m,n\in\mathbb{N}$ are defined by
\[
\Xi_{n,m}=\left\{\begin{array}
[c]{cc}
\Lambda_{n-1}\cdot\Lambda_{n-2}\cdot\cdot\cdot\Lambda_{m} & m<n \\ \mathbf{I} & m=n 
\end{array}
\right .
\]
Furthermore for each $h\leq m\leq n$ we obtain:
\[
 \Xi_{n,m}\Xi_{m,h}=\Xi_{n,h}.
\]
We set with $\mathcal{H}_{\infty}=\underset{\longrightarrow}{\lim}(\mathcal{L}_{n},\ \Xi_{n,m})$ its inductive limit (see \cite{Kadison-Ringrose}) and with 
$Z_{n}:\mathcal{L}_{n}\rightarrow\mathcal{H}_{\infty}$ is the embedding map such that for each natural number $m\leq n$ we have:
\begin{equation}
Z_{n}\Xi_{n,m}=Z_{m}.
\label{Z-map}
\end{equation}
The Hilbert space $\mathcal{H}_{\infty}$ is the closure of linear subspace generated by set $\{Z_{n}\mathcal{L}_{n}:n\in\mathbb{N}\}$, in other words: 
\begin{equation}
\mathcal{H}_{\infty}=\bigvee \limits_{n\in\mathbb{N}}Z_{n}\mathcal{L}_{n}.
\label{Hilbertinfty}
\end{equation}
We observe that the embedding $Z_{n}:\mathcal{L}_{n}\rightarrow\mathcal{H}_{\infty}$ for any $m,n\in\mathbb{N}$ satisfies the following properties:
\[
Z_{n}^{\ast }Z_{m}\left\{ 
\begin{array}{cc}
\Xi _{n,m} & m\leq n \\ 
\Xi _{m,n}^{\ast } & m>n%
\end{array}%
\right. .
\]
We recall that  an isometry dilation of a linear contraction $T$ on Hilbert space $ \mathcal{H} $ (see \cite{Nagy-Foias}) is a  triple $\{\widehat{T},\widehat{ \mathcal{H}},z \}$ with $\widehat{ \mathcal{H}}$ a Hilbert space, $z:\mathcal{H} \rightarrow \widehat{ \mathcal{H}}$  a linear isometry and $\widehat{T}$ an isometry on $\widehat{ \mathcal{H}}$ such that
\[
T^{n}= z^* \widehat{T}^n z \ \ \ \text{for all} \ \  n\in\mathbb{N}.
\]
We observe that if $\Omega$ is a vector belong to $\mathcal{H}$ we have $T\Omega=\Omega$ if and only if $\widehat{T}z\Omega=z\Omega$.
\newline
In fact
\[
||(I-zz^*)\widehat{T}z\Omega||^2= \left\langle(I-zz^*)\widehat{T}z\Omega,(I-zz^*)\widehat{T}z\Omega\right\rangle=
\left\langle \Omega, z^* \widehat{T}^*(I-zz^*)\widehat{T}z\Omega\right\rangle=
\left\langle \Omega, (I-T^*T)\Omega\right\rangle=0.
\]
We can give the following theorem:

\begin{theorem}
\label{CGNS}Let $(\mathfrak{A},\Phi,\varphi)$ be a C*-dynamical system there exist a triple $(\mathcal{H}_{\infty},\pi_{\infty},\Omega_{\infty})$ thus defined:
\begin{itemize}
\item[$\alpha$)] $\mathcal{H}_{\infty}$ is a Hilbert space with $\mathcal{H}_{\varphi}$ embedding in $\mathcal{H}_{\infty}$ i.e. there is a linear isometry $\mathbf{Z}_{0}:\mathcal{H}_{\varphi}\rightarrow\mathcal{H}_{\infty}$;
\item[$\beta$)] $\pi_{\infty}:\mathfrak{A}\rightarrow\mathfrak{B}(\mathcal{H}_{\infty})$ is
a representation such that for each $a\in\mathfrak{A}$ we have
\[
\pi_{\infty}(a)\mathbf{Z}_{0}=\mathbf{Z}_{0}\pi_{\varphi}(a);
\]
\item[$\gamma$)] $\Omega_{\infty}=\mathbf{Z}_{0}\Omega_{\varphi}$.
\end{itemize}
Moreover there exists a linear isometry $\mathbf{V}_{\infty}$ of $\mathfrak{B}(\mathcal{H}_{\infty})$ such that:
\begin{itemize}
\item[1 - ] $\mathbf{V}_{\infty}$ is an isometry dilation of the contraction $\mathbf{U}_{\Phi,\varphi}^*$ :
\[
\mathbf{U}_{\Phi,\varphi}^{n*}=\mathbf{Z}_{0}^*\mathbf{V}_{\infty}^{n}\mathbf{Z}_{0}, \ \ \text{for all} \  n\in\mathbb{N}
\]
and
\[
\mathbf{V}_{\infty}\Omega_{\infty}=\Omega_{\infty};
\]
\item[2 - ] The vector $\Omega_{\infty}$ is cyclic for *-subalgebra $\mathfrak{B}$ of
$\mathfrak{B}(\mathcal{H}_{\infty})$ generated by set 
\[
\bigcup \limits_{n\geq0}\{\mathbf{V}_{\infty}^{n}\pi_{\infty}(a)\mathbf{V}_{\infty}^{* n}:a\in\mathfrak{A}\}.
\]
\item[3 - ] For each $a\in\mathfrak{A}$ we have:
\[
\pi_{\infty}(\Phi(a))=\mathbf{V}_{\infty}^*\pi_{\infty}(a)\mathbf{V}_{\infty}
\]
and
\[
\varphi(a)=\left\langle \Omega_{\infty},\pi_{\infty}(a)\Omega_{\infty}\right\rangle.
\]
\end{itemize}
The quadruple $(\pi_{\infty},\mathcal{H}_{\infty},\Omega_{\infty},\mathbf{V}_{\infty})$  is uniquely determined by the properties 1 - 4 up to unitary equivalence. We shall call any quadruple in this equivalence class, the Covariant GNS representation of the dynamical system $(\mathfrak{A},\Phi,\varphi)$.
\end{theorem}
\begin {proof}
We consider the Hilbert space $\mathcal{H}_{\infty}$ defined in \ref{Hilbertinfty} and with $\mathbf{Z}_{0}:\mathcal{H}_{\varphi}\rightarrow\mathcal{H}_{\infty}$ the linear map \ref{Z-map}.
By the commutative diagram \ref{commutative diagram}, for each natural number
$m\leq n$ and $a$ in $\mathfrak{A}$, we obtain the following relationships:
\begin{itemize}
\item[(e)] $\sigma_{n}(a)\Xi_{n,m}=\Xi_{n,m}\sigma_{m}(a)$;
\item[(f)] $\mathbf{V}_{n}\Xi_{n,m}=\Xi_{n+1,m+1}\mathbf{V}_{m}$;
\item[(g)] $\mathbf{V}_{n}^{\ast}\Xi_{n+1,m}=\Xi_{n,m-1}\mathbf{V}_{m-1}^*$.
\end{itemize}
Then by the properties of inductive limit of directed systems of Hilbert spaces, we can say that there exists a representation $\pi_{\infty}:\mathfrak{A}\rightarrow\mathfrak{B}(\mathcal{H}_{\infty})$ and an isometry $\mathbf{V}_{\infty}:\mathcal{H}_{\infty}\rightarrow\mathcal{H}_{\infty}$ such that:
\begin{itemize}
\item[(h)] $\pi_{\infty}(a)Z_{n}=Z_{n}\sigma_{n}(a)$;
\item[(i)] $\mathbf{V}_{\infty}Z_{n}=Z_{n+1}\mathbf{V}_{n};$
\item[(l)] $\mathbf{V}_{\infty}^*Z_{n}=Z_{n-1}
\mathbf{V}_{n-1}^*$.
\end{itemize}
1) The operator $\mathbf{V}_{\infty}$ is an isometry dilation of the contraction $\mathbf{U}_{\Phi,\varphi}^*$ since for any $h_{n}\in\mathcal{L}_{n}$ we obtain
\[
||\mathbf{V}_{\infty}Z_{n}h_{n}||=||Z_{n+1}\mathbf{V}_{n}h_{n}||=||h_{n}||=||Z_{n}h_{n}||,
\]
while for any $a\in\mathfrak{A}$ and $h\in\mathcal{H}_{\infty}$ we have
\[
\left\langle\pi_{\varphi}(a)\Omega_{\varphi},\mathbf{Z}_{0}^*\mathbf{V}_{\infty}
\mathbf{Z}_{0}h\right\rangle_{\mathcal{H}_{\varphi}}=\left\langle\mathbf{Z}_{0}
\pi_{\varphi}(a)\Omega_{\varphi},\mathbf{Z}_{1}\mathbf{V}_{0}h\right\rangle_{\mathcal{H}_{\infty}} =\left\langle\mathbf{Z}_{1}^*\mathbf{Z}_{0}
\pi_{\varphi}(a)\Omega_{\varphi},\mathbf{V}_{0}h\right\rangle_{\mathcal{L}_{1}} =
\]
\[
=\left\langle \Lambda_{0}\pi_{\varphi}(a)\Omega_{\varphi},\mathbf{V}_{0}h\right\rangle_{\mathcal{L}_{1}} =\left\langle\Omega_{\varphi},\pi_{\varphi}(\Phi(a^*))h\right\rangle_{\mathcal{H}_{\varphi}}=
\left\langle\Omega_{\varphi},\mathbf{U_{\phi,\varphi}^*}h\right\rangle_{\mathcal{H}_{\varphi}}.
\]
Furthermore, for each natural number $n$ we have:
\[
\left\langle\pi_{\varphi}(a)\Omega_{\varphi},\mathbf{Z}_{0}^*\mathbf{V}_{\infty}^{n+1}
\mathbf{Z}_{0}h\right\rangle_{\mathcal{H}_{\varphi}}=  
\left\langle\mathbf{V}_{\infty}^*\mathbf{Z}_{0}\pi_{\varphi}(a)\Omega_{\varphi},\mathbf{V}_{\infty}^{n}
\mathbf{Z}_{0}h\right\rangle_{\mathcal{H}_{\infty}}=
\left\langle\pi_{\varphi}(a)\Omega_{\varphi},\mathbf{U_{\Phi,\varphi}^*}\mathbf{Z}_{0}^*\mathbf{V}_{\infty}^{n}
\mathbf{Z}_{0}h\right\rangle_{\mathcal{H}_{\varphi}},
\]
since $\mathbf{Z}_{0}\pi_{\varphi}(a)\Omega_{\varphi}=\mathbf{Z}_{1}a\overline{\otimes}_{\Phi_{0}}h$, we can write that $\mathbf{V}_{\infty}^*\mathbf{Z}_{0}\pi_{\varphi}(a)\Omega_{\varphi}=
\mathbf{V}_{\infty}^*\mathbf{Z}_{1}a\overline{\otimes}_{\Phi_{0}}h=
\mathbf{Z}_{0}\mathbf{U_{\Phi ,\varphi}}\pi_{\varphi}(a)\Omega_{\varphi}$.
\newline
The vector $\Omega_{\infty}$ is $\mathbf{V}_{\infty}$ invariant, since $\mathbf{U_{\Phi ,\varphi}}^*\Omega_{\varphi}=\Omega_{\varphi}$.
\newline
2) We observe that 
\[
Z_{0}\mathcal{H}_{\varphi}=\overline{\pi_{\infty}((\mathfrak{A})\Omega_{\infty}}
\]
and for each natural number $n$ we obtain
\[
Z_{n}\mathcal{L}_{n}=\overset{n}{\overbrace{\pi_{\infty}(\mathfrak{A})\mathbf{V}_{\infty}\pi_{\infty}(\mathfrak{A})
\mathbf{V}_{\infty} \cdotp\cdotp\cdotp \cdotp \ \pi_{\infty}(\mathfrak{A})\mathbf{V}_{\infty}}}\mathbf{Z}_{0}\mathcal{H}_{\varphi}.
\]
It is also easy to prove that the Hilbert space $Z_{n}\mathcal{L}_{n}$ is generated by follow elements of $\mathcal{H}_{\infty}$:
\[
\partial_{0}(a_{0})\partial_{1}(a_{1}) \cdotp\cdotp\cdotp\cdotp\partial_{n}(a_{n})\Omega_{\infty},
\]
where for each natural number $n$ and $a\in\mathfrak{A}$ we have set
\begin{equation}
\partial_{n}(a)=\mathbf{V}_{\infty}^{n}\pi(a)\mathbf{V}_{\infty}^{n^*}\in\mathfrak{B}.
\end{equation}
3) For each $a\in\mathfrak{A}$ we have
\[\mathbf{V}_{\infty}^*\pi_{\infty}(a)\mathbf{V}_{\infty}Z_{n}=
\mathbf{V}_{\infty}^*\pi_{\infty}(a)Z_{n+1}\mathbf{V}_{n}=\mathbf{V}_{\infty}^*
Z_{n+1}\sigma_{n+1}(a)\mathbf{V}_{n+1}=
\]
\[
=Z_{n}\mathbf{V}_{n}^*\sigma_{n+1}(a)\mathbf{V}_{n}
=Z_{n}\Phi_{n}(a)= Z_{n}\sigma_{n}(\Phi(a))=\pi_{\infty}(\Phi(a))Z_{n},
\]
it follows that
\[
\pi_{\infty}(\Phi(a))=\mathbf{V}_{\infty}^*\pi_{\infty}(a)\mathbf{V}_{\infty} \ \ \text{for all} \ a\in\mathfrak{A}
\]
and
\[
\left\langle\Omega_{\infty},\pi_{\infty}(a)\Omega_{\infty}\right\rangle=
\left\langle \Omega_{\varphi},Z_{0}^*\pi_{\infty}(a)Z_{0}\Omega_{\varphi}\right\rangle=
\left\langle \Omega_{\varphi},\pi_{\varphi}(a)\Omega_{\varphi}\right\rangle=\varphi(a).
\]
Let $(\pi,\mathcal{H},\Omega,\mathbf{V})$ be a new quadruple that satisfies the properties 1 - 4 of the theorem, then there exists an unitary operator $\mathbf{W}:\mathcal{H}_{\infty}\rightarrow \mathcal{H}$ such that
\begin{equation}
\mathbf{W}\pi_{\infty}(a)=\pi(a)\mathbf{W} \label{inductive-limit} \ \ \ \text{for all} \ \ a\in\mathfrak{A}.
\end{equation}
The Hilbert space $\mathcal{L}_{m}$ is generated by elements
\[
\sigma_{m}(a_{m})\mathbf{V}_{m-1}\sigma_{m-1}(a_{m-1}) \cdotp\cdotp\cdotp\cdotp\sigma_{1}(a_{1})\mathbf{V}_{0}\pi_{\varphi}(a_{0})\Omega_{\varphi},
\]
with $a_{1},a_{2},...a_{m}\in\mathfrak{A}$ and  we define a linear operator $\mathbf{W}_{m}:\mathcal{L}_{m}\rightarrow\mathcal{H}$ by
\[
\mathbf{W}_{m}\sigma_{m}(a_{m})\mathbf{V}_{m-1}\sigma_{m-1}(a_{m-1}) \cdotp\cdotp\cdotp\cdotp\sigma_{1}(a_{1})\mathbf{V}_{0}\pi_{\varphi}(a_{0})\Omega_{\varphi}=
\pi(a_{m})\mathbf{V}\pi(a_{m-1})\mathbf{V} \cdotp\cdotp\cdotp\cdotp\pi(a_{1})\mathbf{V}\pi_(a_{0})\Omega,
\]
it is a well defined isometry, since for each $a_{i},b_{j}\in\mathfrak{A}$ with $i,j=1,2,....m$, we obtain
\[
\left\langle\sigma_{m}(a_{m})\mathbf{V}_{m-1}\sigma_{m-1}(a_{m-1}) \cdotp\cdotp\cdotp\cdotp\mathbf{V}_{0}\pi_{\varphi}(a_{0})\Omega_{\varphi},
\sigma_{m}(b_{m})\mathbf{V}_{m-1}\sigma_{m-1}(b_{m-1}) \cdotp\cdotp\cdotp\cdotp\mathbf{V}_{0}\pi_{\varphi}(b_{0})\Omega_{\varphi}
\right\rangle_{\mathcal{L}_{m}}= 
\]
\[=\varphi(a_{0}^*\Phi(a_{1}^*\cdotp\cdotp\cdotp\cdotp\Phi(a_{m-1}^*\Phi(a_{m}^*b_{m})b_{m-1}
\cdotp\cdotp\cdotp\cdotp)b_{1})b_{0})=
\]
\[
=\left\langle\pi(a_{m})\mathbf{V}_{m-1}\pi(a_{m-1})
\cdotp\cdotp\cdotp\mathbf{V}_{0}\pi(a_{0})\Omega,\pi(b_{m})\mathbf{V}_{m-1}\pi(b_{m-1}) \cdotp\cdotp\cdotp\mathbf{V}_{0}\pi(b_{0})\Omega\right\rangle_{\mathcal{H}}.
\]
We observe that for each natural numbers $m,n$ we can write
\[
\Lambda_{n-1}\cdotp\cdotp\cdotp\Lambda_{m}\sigma_{m}(a_{m})\mathbf{V}_{m-1}\sigma_{m-1}(a_{m-1}) \cdotp\cdotp\cdotp\cdotp\sigma_{1}(a_{1})\mathbf{V}_{0}\pi_{\varphi}(a_{0})\Omega_{\varphi}=
\]
\[
=\sigma_{n}(a_{m})\mathbf{V}_{n-1}\sigma_{n-1}(a_{m-1})\cdotp
\cdotp\cdotp\cdotp\sigma_{m+1}(a_{1})\mathbf{V}_{m}\sigma_{m}(a_{0})
\Lambda_{m-1}\cdotp\cdotp\cdotp\cdotp\Lambda_{1}\Lambda_{0}\Omega_{\varphi}.
\]
Then we have the following relationship:
\[
\mathbf{W}_{n}\Xi_{n,m}=\mathbf{W}_{m}
\]
since
\[
\mathbf{W}_{n}\Xi_{n,m}\sigma_{m}(a_{m})\mathbf{V}_{m-1}\sigma_{m-1}(a_{m-1}) \cdotp\cdotp\cdotp\sigma_{1}(a_{1})\mathbf{V}_{0}\pi_{\varphi}(a_{0})\Omega_{\varphi}=
\]
\[
=\mathbf{W}_{n}\sigma_{n}(a_{m})\mathbf{V}_{n-1}\sigma_{n-1}(a_{m-1})\cdotp\cdotp\cdotp\cdotp\cdotp
\sigma_{m+1}(a_{1})\mathbf{V}_{m}\sigma_{m}(a_{0})
\Lambda_{m-1}\cdotp\cdotp\cdotp\Lambda_{1}\Lambda_{0}\Omega_{\varphi},
\]
where 
\[
\Lambda_{m-1}\cdotp\cdotp\cdotp\Lambda_{1}\Lambda_{0}\Omega_{\varphi}=
\mathbf{V}_{m-1}\sigma_{m-1}(1)\mathbf{V}_{m-2}\cdotp\cdotp\cdotp\cdotp\mathbf{V}_{1}\sigma_{1}(1)
\mathbf{V}_{0}\pi_{\varphi}(1)\Omega_{\varphi}.
\]
Therefore
\[
\mathbf{W}_{n}\sigma_{n}(a_{m})\mathbf{V}_{n-1}\sigma_{n-1}(a_{m-1})\cdotp\cdotp\cdotp\cdotp
\sigma_{m+1}(a_{1})\mathbf{V}_{m}\sigma_{m}(a_{0})
\mathbf{V}_{m-1}\sigma_{m-1}(1)\mathbf{V}_{m-2}\cdotp\cdotp\cdotp\mathbf{V}_{1}\sigma_{1}(1)
\mathbf{V}_{0}\pi_{\varphi}(1)\Omega_{\varphi}=
\]
$=\pi(a_{m})\mathbf{V}\pi(a_{m-1}) \cdotp\cdotp\cdotp\cdotp\pi(a_{1})\mathbf{V}\pi_(a_{0})\Omega,\mathbf{V}\pi_(a_{0})\Omega=\mathbf{W}_{m}
\sigma_{m}(a_{m})\mathbf{V}_{m-1}\sigma_{m-1}(a_{m-1}) \cdotp\cdotp\cdotp\sigma_{1}(a_{1})\mathbf{V}_{0}\pi_{\varphi}(a_{0})\Omega_{\varphi}$.
\newline
Moreover we have
\[
\mathcal{H}=\bigvee \limits_{n\in\mathbb{N}}\mathbf{W}_{n}\mathcal{L}_{n},
\]
since $\Omega$ is a cyclic vector for the *-subalgebra $\mathfrak{B}$ of $\mathfrak{B}(\mathcal{H})$ generated by the set $\{\mathbf{V}^n\pi(\mathfrak{A})\mathbf{V}^{* n}:n\in\mathbb{N}\}$.
\newline
Then there is an unitary operator $\mathbf{W}:\mathcal{H}_{\infty}\rightarrow\mathcal{H}$ thus defined:
\[
\mathbf{W}Z_{n}l_{n}=\mathbf{W}_{n}l_{n} \ \ \text{for all} \ l_{n}\in\mathcal{L}_{n},
\]
such that the relationship \ref{inductive-limit} is valid.
\end {proof}

We now turn to some simple observations:
\begin{itemize}
\item [\textit{(a)} -] For each $T$ belong to $\pi_{\infty}(\mathfrak{A})^{\prime}$ we obtain that $z_{o}^*Tz_{o}\in\pi_{\varphi}(\mathfrak{A})^{\prime}$.

\item [\textit{(b)} -] The orthogonal proiection $\mathbf{V}_{\infty}\mathbf{V}_{\infty}^*$ belong to $\pi_{\infty}(\mathcal{D}_{\Phi})^{\prime}$, 
where $\mathcal{D}_{\Phi}$ is the multiplicative domain of the ucp-map $\Phi$ and
\[
\mathbf{V}_{\infty}^*\pi_{\infty}(a)=\pi_{\infty}(\Phi(a))\mathbf{V}_{\infty}^* \ \ \text{for all} \ a\in\mathcal{D}_{\Phi}.
\]
In fact for each natural number $n$, we have that  $\mathcal{D}_{\Phi}\subset\mathcal{D}_{\Phi_{n}}$ with $\mathcal{D}_{\Phi_{n}}$ the multiplicative domains of the ucp-map $\Phi_{n}:\mathfrak{A}\rightarrow\mathfrak{B}(\mathcal{L}_{n})$ previous defined in \ref{ucpn-map}, hence
\[
\mathcal{D}_{\Phi}=\bigcap\limits_{n\geq0}\mathcal{D}_{\Phi_{n}}.
\]
\item [\textit{(c)} -] If $\Phi$ is an automorphism, the triple $\left(  \mathcal{H}_{\infty}%
,\pi_{\infty},\Omega_{\infty}\right)  $ is a unitary equivalent to the
GNS\ $\left(  \mathcal{H}_{\varphi},\pi_{\varphi},\Omega_{\varphi}\right)  $.

\item [\textit{(d)} -] If $\varphi$ is a faithful state of $\mathfrak{A}$ we have that
$\pi_{\infty}:\mathfrak{A}\rightarrow\mathfrak{B}(\mathcal{H}_{\infty})  $ is a faithful representation and $\Omega_{\infty}$ is a separating vector for $\pi_{\infty}(\mathfrak{A})$.

\item [\textit{(e)} -] If $(\mathfrak{M},\Phi,\varphi)$ is a W*-dynamical system, the CGNS representation
$\pi_{\infty}:\mathfrak{M}\rightarrow\mathfrak{B}(\mathcal{H}_{\infty})$ 
is faithful and normal, since the Stinespring representations 
$(\mathcal{L}_{n+1},\sigma_{n+1},\mathbf{V}_{n})$ previous defined, are normal maps for all $n\in\mathbb{N}$.
\end{itemize}

We now study the covariant GNS representation for $C^*$-dynamical systems with multiplicative dynamics (see proposition 6.2 in \cite{NSZ}).
\begin{proposition}
Let $(\mathfrak{A},\Phi,\varphi)  $ be a C*-dynamical system with $\Phi$ homomorphism and
$(\pi_{\infty},\mathcal{H}_{\infty},\Omega_{\infty},\mathbf{V}_{\infty
})  $ the CGNS representation described in Theorem \ref{CGNS}, we obtain that $\mathbf{V}_{\infty}:\mathcal{H}_{\infty
}\rightarrow\mathcal{H}_{\infty}$ is a unitary operator, since for each
$n$ in $\mathbb{N}$ the linear isometry 
$\mathbf{V}_{n}:\mathcal{L}_{n}\rightarrow\mathcal{L}_{n+1}$ 
of the Stinespring representations $(\mathcal{L}_{n+1},\sigma_{n+1},\mathbf{V}_{n})  $ are unitary operators. We can write $\mathbf{V}_{n}\mathcal{L}_{n}=\mathcal{L}_{n+1}$ with
\[
\mathbf{V}^nZ_{0}\mathcal{H}_{\varphi}=Z_{n}\mathcal{L}_{n}.
\]
Moreoverer $\mathbf{V}_{\infty}$ is the minimal unitary dilation of $\mathbf{U}_{\Phi,\varphi}^*$:
\begin{equation}
\mathcal{H}_{\infty}=\bigvee \limits_{n\in\mathbb{N}}
\mathbf{V}_{\infty}^{n}\mathbf{Z}_{0}\mathcal{H}_{\varphi},
\label{prop-minimal2}
\end{equation}
and 
\begin{equation}
\pi_{\infty}(a)\mathbf{V}_{\infty}=\mathbf{V}_{\infty}\pi_{\infty}(\Phi(a)) \ \ \text{for all} \ a\in\mathfrak{A}.
\end{equation}
\end{proposition}
\begin{proof}
it's a trivial consequences of the previous propositions.
\end{proof}
Finally, for C*-dynamical system with multiplicative dynamics  we gives the following result:
\begin{lemma}
Let $(\mathfrak{A},\Phi,\varphi)$ be a C*-dynamical system with $\Phi$ a homomorphism and $\Omega_{\varphi}$ cyclic vector for 
$\pi_{\varphi}(\mathfrak{A})^{\prime}$, its CGNS representation $(\pi_{\infty},\mathcal{H}_{\infty},\Omega_{\infty},\mathbf{V}_{\infty
})$ has the following properties:
\begin{enumerate} 
\item[1 - ] $\pi_{\varphi}(a)=0$ if and only if $\pi_{\infty}(a)=0$;
\item[2 - ]  $\Omega_{\infty}$ is a separanting vector for $\pi_{\infty}(\mathfrak{A})$;
\item[3 - ] $|| \pi_{\infty}(a)||=|| \pi_{\varphi}(a)||$ \ \ \ \ for all $a\in\mathfrak{A}$.

\end{enumerate}
\end{lemma}
\begin {proof}
1) For each natural number $n$ we have that $\pi_{\infty}(a)\mathbf{V}_{\infty}^nz_{0}=0$ and from the relationship \ref{prop-minimal2} follows that the $\pi_{\infty}(a)=0$.
\newline
In fact:
\[
\pi_{\infty}(a)\mathbf{V}_{\infty}^nz_{0}\Omega_{\varphi}=
\pi_{\infty}(a)\mathbf{V}_{\infty}^n\Omega_{\infty}=
\pi_{\infty}(a)\Omega_{\infty}=
Z_{0}\pi _{\varphi}(a)\Omega_{\varphi}=0.
\]
Moreover 
\[
\pi_{\varphi}(\Phi^n(a^*a))\Omega_{\varphi}=
z_{0}^*\mathbf{V}_{\infty}^{n *}\pi_{\infty}(a^*a)\mathbf{V}_{\infty}^{n}z_{0}\Omega_{\varphi}=0,
\]
with $\Omega_{\varphi}$ is a separating vector for the von Neumann algebra $\pi_{\varphi}(\mathfrak{A})^{\prime \prime}$, hence we obtain 
\[
z_{0}^*\mathbf{V}_{\infty}^{n *}\pi_{\infty}(a^*)\pi_{\infty}(a)\mathbf{V}_{\infty}^{n}z_{0}=0.
\]
2) If $\pi_{\infty}(a)\Omega_{\infty}=0$ it follows that 
\[
\pi_{\infty}(a)\Omega_{\infty}=\pi_{\infty}(a)Z_{0}\Omega_{\varphi}=Z_{0}\pi_{\varphi}(a)\Omega_{\varphi}=0,
\]
since $Z_{0}:\mathcal{H}_{\varphi}\rightarrow\mathcal{H}_{\infty}$ is an isometric operator and $\Omega_{\varphi}$ is a separating vector for $\pi_{\varphi}(\mathfrak{A})^{\prime \prime}$, we obtain that $\pi_{\varphi}(x)=0$, then $\pi_{\infty}(x)=0$.
\newline
3) Obviously, for each $a$ belong to $\mathfrak{A}$  we have $||\pi_{\varphi}(a)||\leq||\pi_{\infty}(a)||$ since $\pi_{\varphi}(a)=z_{0}^*\pi_{\infty}(a)z_{o}$.
\newline
By the second statement of the proposition, $\Omega_{\infty}$ is a separating vector for the von Neumann algebra  $\pi_{\infty}(\mathfrak{A})^{\prime}$ and for each $a$ belong to $\mathfrak{A}$ and $T$ in $\pi_{\infty}(\mathfrak{A})^{\prime}$ we can write:
\[
||\pi_{\infty}(a)T\Omega_{\infty}||\leq||T\Omega_{\infty}|| ||\pi_{\varphi}(a)||,
\]
hence we obtain that $||\pi_{\infty}(a)||\leq||\pi_{\varphi}(a)||$.
\newline
In fact 
\[
||\pi_{\infty}(a)T\Omega_{\infty}||^2=\left\langle\Omega_{\varphi},z_{o}^*T^*\pi_{\infty}(a^*)\pi_{\infty}(a)
Tz_{o}\Omega_{\varphi}\right \rangle=
\left\langle\Omega_{\varphi},z_{o}^*T^*Tz_{o}\pi_{\varphi}(a^*a)\Omega_{\varphi}\right \rangle.
\]
The positive element $z_{o}^*T^*Tz_{o}$ belong to von Neumann algebra $\pi_{\varphi}(\mathfrak{A})^{\prime}$ it follows that there is a element $Y$ in $\pi_{\varphi}(\mathfrak{A})^{\prime}$ such that
$z_{o}^*T^*Tz_{o}=Y^*Y$.
\newline
Then we obtain:
\[
||\pi_{\infty}(a)T\Omega_{\infty}||^2=
\left\langle\Omega_{\varphi},Y^*Y\pi_{\varphi}(a^*a)\Omega_{\varphi}\right \rangle=
\left\langle\Omega_{\varphi},Y^*\pi_{\varphi}(a^*)Y\pi_{\varphi}(a)\Omega_{\varphi}\right \rangle=
\]
\[
=||\pi_{\varphi}(a)Y\Omega_{\varphi}||^2
 \leq||Y\Omega_{\varphi}||^2||\pi_{\varphi}(a)||^2,
\]
with $||Y\Omega_{\varphi}||^2=||T\Omega_{\infty}||^2$.

\end {proof}

\section{Reversible dilation for C*-dynamical systems with multiplicative dynamics} \label{REVERSIBLE-DIL}
In this section we will use the CGNS representation to prove that the W*-dynamical system associated to C*-dynamical system with multiplicative dynamics, admits a minimal reversible dilation, that keeps unchanged the ergodic properties of the original system. Furthermore we shall show that a dynamical system which admits a right inverse, admit a minimal reversible dilation.
\newline
Let us briefly summarize the main concepts and results needed in this section.
\newline
Let $(\mathfrak{A},\Phi,\varphi)$ be a C*-dynamical system, we say that the ucp-map $\Phi$ admit a $\varphi$-adjoint, if
there is an ucp-map $\Phi^{\sharp}:\mathfrak{A}\rightarrow\mathfrak{A}$ such
that for any $a,b\in\mathfrak{A}$
\[
\varphi(a\Phi^{\sharp}(b))=\varphi(\Phi(a)b).
\]
The property of adjunction of a state, fundamental in reversible processes and ergodic theory, has been studied by various authors (See e.g.\ \cite{Antha}, \cite{Majewski} and \cite{NSZ}) and its basic properties are summarized in the following proposition:
\begin{proposition}
\label{prop-NSZ}
Let $(\mathfrak{A},\Phi,\varphi)$ be a C*-dynamical system with $\Omega_{\varphi}$ cyclic for the von Neumann algebra $\pi_{\varphi}(\mathfrak{A})^{\prime}$ and  $(\Delta_{\varphi},\mathbf{J}_{\varphi})$ the modular operators associated with pair $(\pi_{\varphi}(\mathfrak{A})^{\prime\prime},\Omega_{\varphi})$. 
\newline
The following conditions are equivalent:
\begin{itemize}
\item[1 - ] $\Phi$ commutes with the automorphism modular group i.e. 
\[
\sigma_{t}^{\varphi}\circ\Phi_{\bullet}=\Phi_{\bullet}\circ\sigma_{t}^{\varphi} \ \ 
\text{for all} \ t\in\mathbb{R};
\]
\item[2 - ] $\mathbf{U}_{\Phi,\varphi}$ commutes with modular operators:
\[
\mathbf{U}_{\Phi,\varphi}\Delta_{\varphi}^{it}=\Delta_{\varphi}^{it}\mathbf{U}_{\Phi,\varphi}  \ \ \text{for all} \ t\in\mathbb{R};
\]
and
\[
\mathbf{U}_{\Phi,\varphi}\mathbf{J}_{\varphi}=\mathbf{J}_{\varphi
}\mathbf{U}_{\Phi,\varphi}; 
\]
\item[3 - ] There exists an unique normal ucp-map $\Phi^{\sharp}:\pi_{\varphi}(\mathfrak{A})^{\prime\prime}\rightarrow\pi_{\varphi}(\mathfrak{A})
^{\prime\prime}$ such that for each $a\in\mathfrak{A}$ we have
\[
\mathbf{U}_{\Phi,\varphi}^{\ast}\pi_{\varphi}(a)\Omega_{\varphi}=\pi_{\varphi}(\Phi^{\sharp}(a))\Omega_{\varphi}.
\]
\item[4 - ] If the dynamics $\Phi$ is a homomorphism, then the previous conditions are equivalent also with the following:
\[
\mathbf{U}_{\Phi,\varphi}^*\pi_{\varphi }(\mathfrak{A})^{\prime \prime}\mathbf{U}_{\Phi,\varphi }
\subset\pi_{\varphi}(\mathfrak{A})^{\prime\prime}.
\]
\end{itemize}
\end{proposition}
\begin {proof}
See proposition 3.3 in \cite{NSZ}.
\end {proof}
We give now the definition of reversible dilation of a W*-dynamical system (see \cite{Kummerer}):
\begin{definition}
A W*-dynamical system $(\widehat{\mathfrak{M}},\widehat{\Phi},\widehat{\varphi})$ with dynamics $\widehat{\Phi}$ an automorphism, is said to be a reversible dilation of the W*-dynamical system $(\mathfrak{M},\Phi,\varphi)$, if it satisfies the following conditions:
\newline
There is a normal ucp-map $\mathcal{E}:\widehat{\mathfrak{M}}\rightarrow \mathfrak{M}$ and a normal injective homomorphism $i:\mathfrak{M}\rightarrow \widehat{\mathfrak{M}}$ such that for each $a$ belong to $\mathfrak{A}$ and $X$ in $\widehat{\mathfrak{M}}$ we have:
\[
\mathcal{E}(i(a)X)=a\mathcal{E}(X)
\]
and for each natural number $n$
\[
\mathcal{E}(\widehat{\Phi }^{n}(i(a)))=\Phi ^{n}((a)),
\]
with
\[
\widehat{\varphi}(X)=\varphi(\mathcal{E}(X)).
\]
Furthermore the dilation is said to be minimal if the von Neumann algebra $\widehat{\mathfrak{M}}$ is generated by the set:
\[
\bigcup\limits_{k\in \mathbb{Z}}\{\widehat{\Phi }^{k}(i(a)):a\in\mathfrak{M}\}.
\]
\end{definition}
We observe that the ucp-map $\widehat{\mathcal{E}} = i\circ \mathcal{E}$ is a conditional expectation from  $\widehat{\mathfrak{M}}$ onto $i(\mathfrak{M})$ which leave invariant a faithful normal state. The existence of such map be derived from a theorem of Takesaki (see \cite{Take} and for its generalization \cite{Accardi-Cecchini}) which characterize the range of existence of a reversible dilation of a dynamical system.
\newline
Furthermore, it is easily to show  that if a W*-dynamic system $(\mathfrak{M},\Phi,\varphi)$ admit a reversible dilation, the dynamic $\Phi$ admit a $\varphi $-adjoint $\Phi^{\sharp}$ (see \cite{Kummerer}).
\newline
 
The following result is a reformulation of proposition 6.2 in \cite{NSZ}, we include a proof for completeness.
\begin{proposition}
Let $(\mathfrak{A},\varphi,\Phi)$ be a C*-dynamical system with $\Phi$ homomorphism and $(\mathcal{H}_{\infty},\pi_{\infty},\Omega_{\infty},\mathbf{V}_{\infty })$ its CGNS representation. If $\Omega_{\varphi }$ is a cyclic vector for $\pi_{\varphi}( \mathfrak{A})^{\prime}$ we have that $\Omega_{\infty }$ is a separating vector for the von Neumann algebra 
$\mathfrak{B}^{^{\prime \prime }}$, where $\mathfrak{B}$ is the unital *-subalgebra of $\mathfrak{B}(\mathcal{H}_{\infty })$ generated by the set:
\begin{equation}
\bigcup\limits_{n\in \mathbb{N}}\{\mathbf{V}_{\infty }^{n}\pi_{\infty}(a)
\mathbf{V}_{\infty }^{* n}: a\in\mathfrak{A}\}. \label{algebra-B}
\end{equation}
Moreover if $\Phi$ admit a $\varphi$-adjoint $\Phi^{\sharp}$, there is a normal ucp-map  $\mathcal{E}:\mathfrak{B}^{\prime \prime}\rightarrow\pi_{\varphi}(\mathfrak{A})^{\prime \prime}$ such that for each natural number $n$ and element $a\in\mathfrak{A}$ we have
\begin{equation}
\mathcal{E}(\mathbf{V}_{\infty}^{n}\pi_{\infty}(a)\mathbf{V}_{\infty}^{n *})=\pi_{\varphi}
(\Phi^{\sharp n}(a))   \label{econd.1}.
\end{equation}
Furthermore
\begin{equation}
\mathcal{E}(\pi_{\infty}(a)X)=\pi_{\varphi}(a)\mathcal{E}(X) \label{econd.2}
\end{equation}
and for any $a\in\mathfrak{A}$ and $X\in\mathfrak{B}^{\prime \prime}$ we obtain
\begin{equation}
\left\langle \Omega_{\varphi}, \mathcal{E}(X) \Omega_{\varphi} \right \rangle=
\left\langle \Omega_{\infty},X \Omega_{\infty} \right \rangle. \label {econd.3}
\end{equation}
\end{proposition}

\begin {proof}
We observe that for each natural number $n$ we obtain the following inclusions:
\[
\pi_{\infty}(\mathfrak{A})\subset\mathbf{V}_{\infty}\pi_{\infty}(\mathfrak{A})\mathbf{V}_{\infty}^*
\subset\mathbf{V}_{\infty}^{2}\pi_{\infty}(\mathfrak{A})\mathbf{V}_{\infty}^{*2}\subset
\cdot\cdot\cdot\subset\mathbf{V}_{\infty}^{n}\pi_{\infty}(\mathfrak{A})\mathbf{V}_{\infty}^{*n}\subset\cdot\cdot\cdot.
\]
In fact
$\mathbf{V}_{\infty}^*\pi_{\infty}(\mathfrak{A})\mathbf{V}_{\infty}=\pi_{\infty}(\Phi(\mathfrak{A}))  \subset\pi_{\infty}(\mathfrak{A})$,
with $\mathbf{V}_{\infty}$ an unitary operator, so we can write that
$\pi_{\infty}(\mathfrak{A})\subset\mathbf{V}_{\infty}\pi_{\infty}(\mathfrak{A})\mathbf{V}_{\infty}^*$.
\newline
Then, let $X$ be any element belong to the *-algebra $\mathfrak{B}$, we can write it as follows:
\[
X=\mathbf{V}_{\infty}^{k}\pi_{\infty}(x)\mathbf{V}_{\infty}^{* k}
\]
for some $x\in\mathfrak{A}$ and $k\in\mathbb{N}$.
\newline
We observe that for each natural number $k$ and $x\in\mathfrak{A}$ we have:
\[
Z_{0}^*\mathbf{V}^{k}\pi(x)\mathbf{V}^{k *}Z_{0}=\mathbf{U}_{\Phi,\varphi }^{k *}\pi _{\varphi}(x)\mathbf{U}_{\Phi ,\varphi}^k
\]
and from the proposition \ref{prop-NSZ}, we can say that
 $Z_{0}^*XZ_{0}\in\pi _{\varphi }(\mathfrak{A})^{\prime\prime}$ for all $X\in\mathfrak{B}$.
\newline
Then $Z_{o}^*XZ_{0}\Omega_{\varphi}=\mathbf{U}_{\Phi,\varphi }^{k *}\pi _{\varphi}(x)\mathbf{U}_{\Phi ,\varphi}^k\Omega_{\varphi}=\pi_{\varphi}(\Phi^{\sharp k}(x))\Omega_{\varphi}$ 
with $\Omega_{\varphi}$ separating vector for $\pi_{\varphi }(\mathfrak{A})^{\prime\prime}$, hence we obtain
\[
Z_{0}^*\mathbf{V}^{k}\pi(x)\mathbf{V}^{k *}Z_{0}=\pi_{\varphi}(\Phi^{\sharp k}(x)).
\]
Furthermore, by the double commutant theorem 
\[
Z_{0}^*\mathfrak{B}^{\prime \prime}Z_{0}\subset\pi_{\varphi}(\mathfrak{A})^{\prime \prime},
\]
therefore we obtain a normal ucp-map $\mathcal{E}:\mathfrak{B}^{\prime \prime}\rightarrow\pi_{\varphi}(\mathfrak{A})^{\prime \prime}$ such that
\[ 
\mathcal{E}(X)=Z_{0}^*XZ_{0} \ \ \text{for all} \ \ X\in\mathfrak{B}^{\prime \prime}
\]
with $\mathcal{E}(\mathfrak{B})\subset\pi_{\varphi}(\mathfrak{A})$.
\newline
Moreover, for each $a\in\mathfrak{A}$ and $X\in\mathfrak{B}^{\prime \prime}$ we have
\[
\mathcal{E}(\pi_{\infty}(a)X)=Z_{0}^*\pi_{\infty}(a)XZ_{0}=\pi_{\varphi}(a)Z_{0}^*XZ_{0}=
\pi_{\varphi}(a)\mathcal{E}(X).
\]
We observe that if $X$ belong to $\mathfrak{B}$ we have $\mathbf{V}^{* n}_{\infty}X\mathbf{V}_{\infty}^n\in\mathfrak{B}$ for all $n\in\mathbb{N}$, therefore from double commutant theorem, for each natural number $n$ we can write
\[
\mathbf{V}^{* n}_{\infty}\mathfrak{B}^{\prime \prime}\mathbf{V}_{\infty}^n\subset\mathfrak{B}^{\prime \prime}.
\]
If $X$ belong to $\mathfrak{B}^{\prime \prime}$ with $X\Omega_{\infty}=0$ we have
\[
X\Omega_{\infty}=X\mathbf{V}_{\infty}^{n}\Omega_{\infty}=X\mathbf{V}_{\infty}^{n}Z_{0}\Omega_{\varphi}=0
\]
for all  $n\subset\mathbb{N}$. It follows that 
$Z_{0}^*\mathbf{V}_{\infty}^{* n}X^*X\mathbf{V}_{\infty}^nZ_{0}\Omega_{\varphi}=0$
with $\Omega_{\varphi}$ separable vector for von Neumann algebra $\pi_{\varphi}(\mathfrak{A})^{\prime \prime}$ and $\mathbf{V}_{\infty}^{* n}X^*X\mathbf{V}_{\infty}^n\in\mathfrak{B}^{\prime \prime}$, hence $X\mathbf{V}_{\infty}^nZ_{0}=0$ and from the relationship \ref{prop-minimal2} we obtain that $X=0$.
\end {proof}
Our main tool in this section is the following proposition
\begin{proposition}
Let $(\mathfrak{A},\Phi,\varphi )$ be a C*-dynamical system with $\Phi$ homomorphism and $\Omega_{\varphi }$ cyclic vector for $\pi_{\varphi}(\mathfrak{A})^{\prime}$. 
If $\Phi$ admit a $\varphi$-adjoint $\Phi^{\sharp}$, the W*-dynamical system $(\pi_{\varphi}(\mathfrak{A})^{\prime\prime},\Phi_{\bullet},\varphi_{\bullet})$ associated to our C* dynamical system, admit a minimal reversible dilation  $(\widehat{\mathfrak{M}},\widehat{\Phi},\widehat{\varphi},i,\mathcal{E})$ where:
\begin{itemize}
\item[1  -] 
The von Neumann algebra $\widehat{\mathfrak{M}}$ is double commutant of the*-subalgebra $\mathfrak{B}$ previous defined in \ref{algebra-B};
\item[2 -] The injective homomorphism $i$ is thus defined:
\[
i(A)\Omega_{\infty}=z_{o}A\Omega_{\varphi} \ \ \text{for all} \  A\in\pi_{\varphi}(\mathfrak{A})^{\prime \prime}, 
\]
while for the automorphism $\widehat{\Phi}:\widehat{\mathfrak{M}}\rightarrow\widehat
{\mathfrak{M}}$ we have:
\[
\widehat{\Phi}(X)=\mathbf{V}_{\infty}^*X\mathbf{V}_{\infty} \ \ \text{for all} \ X\in\widehat{\mathfrak{M}};
\]

\item[3 -]
The conditional expectation $\mathcal{E}:\widehat{\mathfrak{M}}\rightarrow\mathfrak{M}$ is defined through the relationship:
\[
\mathcal{E}(\widehat{\mathbf{V}}_{\infty}^{k}\pi_{\infty}(a)\widehat{\mathbf{V}}_{\infty}^{-k})=
\Phi^{\sharp k}(a)
\]
for all $a\in\mathfrak{A}$ and $k\in\mathbb{N}$, while for the faithful normal state $\widehat{\varphi}$ we have:
\[
\widehat{\varphi}(X) = \left\langle \Omega_{\infty}, X \Omega_{\infty} \right \rangle \ \ \text{for all} \  X\in\widehat{\mathfrak{M}}.
\]
\end{itemize}
\end{proposition}
\begin {proof}
The isometric homomorphism $i_{o}:\pi_{\varphi}(\mathfrak{A})\rightarrow\pi_{\infty}(\mathfrak{A})$ thus defined 
\[
i_{o}(\pi_{\varphi}(x))=\pi_{\infty}(x) \ \ \text{for all} \ x\in\mathfrak{A}, 
\]
can be uniquely extended to a normal homomorphism 
$i:\pi_{\varphi}(\mathfrak{A})^{\prime \prime}\rightarrow\pi_{\infty}(\mathfrak{A})^{\prime \prime}$.
\newline
In fact, let $A\in\pi_{\varphi}(\mathfrak{A})^{\prime \prime}$, by the Kaplansky density theorem there is a net $\{a_{\alpha}\}_{\alpha}$ satisfying 
$||a_{\alpha}|| \leq ||A||$ and $\pi_{\varphi}(a_{\alpha})\rightarrow A$ in $so-top$ (\textit{i.e.} in strong operator topology), hence we obtain $ z_{o}\pi_{\varphi}(a_{\alpha})\Omega_{\varphi}\rightarrow  z_{o}A\Omega_{\varphi}$.
\newline
The closed unit ball of von Neumann algebra is $\omega-top$ compact (with $\omega-top$ we set normal topology of a von Neumann algebra), let $X$ be any $\omega-top$ limit point of the bounded net $\{\pi_{\infty}(a_{\alpha})\}_{\alpha}$ we obtain that $X\Omega_{\infty}=z_{o}A\Omega_{\varphi}$  since $\pi_{\infty}(a_{\alpha})\Omega_{\infty}=z_{o}\pi_{\varphi}(a_{\alpha})\Omega_{\varphi}$, then $X$ is a unique $\omega-top$ limit point and we can define:
\[
i(A)=X.
\]
We need only to prove that the diagram 
\[
\begin{array}{ccccc}
\widehat{\mathfrak{M}} &  & \overset{\widehat{\Phi }^{n}}{\longrightarrow }
&  & \widehat{\mathfrak{M}} \\ 
& 
\begin{array}{c}
\ \widehat{\varphi } \\ 
\searrow 
\end{array}
&  & 
\begin{array}{c}
\widehat{\varphi }\  \\ 
\swarrow 
\end{array}
&  \\ 
i \uparrow  &  & \mathbb{C} &  & \downarrow \mathcal{E} \\ 
& 
\begin{array}{c}
\varphi  \\ 
\nearrow 
\end{array}
&  & 
\begin{array}{c}
\varphi  \\ 
\nwarrow 
\end{array}
&  \\ 
\pi_{\varphi}(\mathfrak{A})^{\prime \prime} &  & \overset{\Phi_{\bullet}^{n}}{\longrightarrow } &  & \pi_{\varphi}(\mathfrak{A})^{\prime \prime}
\end{array}%
\]
is commutative for all natural number $n$.
\newline
In fact, for each $A$ belong to $\pi_{\varphi}(\mathfrak{A})^{\prime \prime}$ we have:
\[
\mathcal{E}(\widehat{\Phi}^n(i(A))\Omega_{\varphi}=
z_{o}^*\mathbf{V}_{\infty}^{n *}i(A)\mathbf{V}_{\infty}^nz_{o}\Omega_{\varphi}=
z_{o}^*\mathbf{V}_{\infty}^{n *}i(A)\Omega_{\infty}=
z_{o}^*\mathbf{V}_{\infty}^{n *}z_{o}A\Omega_{\varphi}=
\mathbf{U}_{\Phi, \varphi}^{n}A\Omega_{\varphi}=\Phi_{\bullet}^n(A)\Omega_{\varphi}
\]
while for each $X\in\mathfrak{B}^{\prime \prime}$ we obtain:
\[
\varphi_{\bullet}(\mathcal{E}(X))=
\left\langle\Omega_{\varphi}, z_{o}^*Xz_{o}\Omega_{\varphi}\right\rangle=
\left\langle\Omega_{\infty}, X\Omega_{\infty}\right\rangle=\widehat{\varphi}(X)
\]
and
\[
\widehat{\varphi}(\widehat{\Phi}(X))=
\left\langle\Omega_{\infty}, \mathbf{V}_{\infty}^*X\mathbf{V}_{\infty}\Omega_{\infty}\right\rangle=
\widehat{\varphi}(X).
\]
\end {proof}
In finally we have the following remark:
\begin{remark}
Any W*-dynamic system $(\mathfrak{M},\Phi,\varphi )$ with dynamics $\Phi$ a homomorphism, admit a reversible dilation.
\end{remark}

We study now the ergodic properties of the dilation, determined by the previous proposition.
\newline
We recall that a C*-dynamical system $(\mathfrak{A},\Phi,\varphi)$ is ergodic if
\begin{equation}
\underset{n\rightarrow\infty}{\lim}\dfrac{1}{n+1}
\sum \limits_{k=0}^{n}[\varphi(a\Phi^{k}(b))-\varphi(a)\varphi(b)]=0, \ \ \text{for all} \ a,b\in\mathfrak{A},
\end{equation}
 while it is a weakly mixing if
\begin{equation}
\underset{n\rightarrow\infty}{\lim}\dfrac{1}{n+1}
\sum \limits_{k=0}^{n}|\varphi(a\Phi^{k}(b))-\varphi(a)\varphi(b)|=0, \ \text{for all} \ a,b\in\mathfrak{A}.
\end{equation}
A most general and abstract framework for the study of the noncommutative ergodic theory is found in \cite{NSZ}.
\begin{proposition}
Let $(\mathfrak{A},\Phi,\varphi )$ be a C*-dynamical system with dynamics $\Phi$ a homomorphism that admit a $\varphi$-adjoint and $\Omega_{\varphi }$ 
cyclic vector for $\pi_{\varphi}(\mathfrak{A})^{\prime}$. If the C*-dynamical system is ergodic [weakly mixing], the reversible dilation of the previous proposition, of its associated W*-dynamical system $(\pi_{\varphi}(\mathfrak{A})^{\prime\prime},\Phi_{\bullet},\varphi_{\bullet})$ is ergodic [weakly mixing].
\end{proposition}
\begin {proof}
We will prove that for each $X,Y\in\mathfrak{B}^{\prime \prime}$ result
\[
\lim_{N\rightarrow\infty}\frac{1}{N+1}\sum \limits_{k=0}^{N}
\left[\widehat{\varphi}(X\widehat{\Phi}^{k}(Y))-\widehat{\varphi}(X)\widehat{\varphi}(Y)\right]=0 .
\]
Let $X\in\mathfrak{B}^{\prime \prime}$ and $Y\in\mathfrak{B}$ with $Y=\mathbf{V}_{\infty}^{j}\pi_{\infty}(y)\mathbf{V}_{\infty}^{j *}$, we have for any  $k\geq j$:
\[
\widehat{\varphi}(X\widehat{\Phi}^{k}(Y))=
\left\langle\Omega_{\varphi}, z_{o}^*X\mathbf{V}_{\infty}^{(k-j)*}\pi_{\infty}(y)z_{o}\Omega_{\varphi}\right\rangle=
\left\langle\Omega_{\varphi}, z_{o}^*X\mathbf{V}_{\infty}^{(k-j)*}z_{o}\pi_{\varphi}(y)\Omega_{\varphi}\right\rangle=
\]
\[
=\left\langle\Omega_{\varphi},
z_{o}^*Xz_{o}\mathbf{U}_{\Phi \varphi}^{(k-j)*}\pi_{\varphi}(y)\Omega_{\varphi}\right\rangle= 
\left\langle\Omega_{\varphi},
z_{o}^*Xz_{o}\pi_{\varphi}(\Phi^{(k-j)}(y))\Omega_{\varphi}\right\rangle= 
\left\langle\Omega_{\varphi},
\mathcal{E}(X)\Phi_{\bullet}^{(k-j)}(y)\Omega_{\varphi}\right\rangle.
\]
It follows that
\[
\lim_{N\rightarrow\infty}\frac{1}{N+1}\sum \limits_{k=0}^{N}
[\widehat{\varphi}(X\widehat{\Phi}^{k}(Y))-\widehat{\varphi}(X)\widehat{\varphi}(Y)] =
\lim_{N\rightarrow\infty}\frac{1}{N+1}\sum \limits_{k=0}^{N}
[\varphi_{\bullet}(\mathcal{E}(X)\Phi_{\bullet}^{(k-j)}(y))-
\varphi_{\bullet}(\mathcal{E}(X))\varphi_{\bullet}(\mathcal{E}(Y)] =
\]
\[
=\lim_{N\rightarrow\infty}\frac{1}{N+1}\sum \limits_{k=0}^{N}
[\varphi_{\bullet}(\mathcal{E}(X)\Phi_{\bullet}^{k}(y))-
\varphi_{\bullet}(\mathcal{E}(X))\varphi_{\bullet}(\mathcal{E}(Y)] =0.
\]
Let $Y\in\mathfrak{B}^{\prime \prime }$, for each $\epsilon>0$ there is a element $Y_{\epsilon}\in\mathfrak{B}$ such that $||Y-Y_{\epsilon}||< \epsilon$.
\newline
Then 
\[
\widehat{\varphi}(X\widehat{\Phi}^{k}(Y))-\widehat{\varphi}(X)\widehat{\varphi}(Y)= \widehat{\varphi}(X\widehat{\Phi}^{k}(Y_{\epsilon}))-\widehat{\varphi}(X)\widehat{\varphi}(Y_{\epsilon})+
\widehat{\varphi}(X\widehat{\Phi}^{k}(Y-Y_{\epsilon}))-\widehat{\varphi}(X)\widehat{\varphi}(Y-Y_{\epsilon})
\]
and 
\[
\lim_{N\rightarrow\infty}\frac{1}{N+1}\sum \limits_{k=0}^{N}
[\widehat{\varphi}(X\widehat{\Phi}^{k}(Y))-\widehat{\varphi}(X)\widehat{\varphi}(Y)] =
\lim_{N\rightarrow\infty}\frac{1}{N+1}\sum \limits_{k=0}^{N}
[\widehat{\varphi}(X\widehat{\Phi}^{k}(Y_{\epsilon}))-\widehat{\varphi}(X)\widehat{\varphi}(Y_{\epsilon})]+
\]
\[
+\lim_{N\rightarrow\infty}\frac{1}{N+1}\sum \limits_{k=0}^{N}
[\widehat{\varphi}(X\widehat{\Phi}^{k}(Y-Y_{\epsilon}))-\widehat{\varphi}(X)\widehat{\varphi}(Y-Y_{\epsilon})]=0.
\]
Since $||\widehat{\Phi}^k||=1$ for all natural number $k$ and 
\[
\frac{1}{N+1}|\sum \limits_{k=0}^{N}
[\widehat{\varphi}(X\widehat{\Phi}^{k}(Y-Y_{\epsilon}))-\widehat{\varphi}(X)\widehat{\varphi}(Y-Y_{\epsilon})]|
\leq 2\epsilon ||X||.
\]
The proof of the weakly mixing is performed in the same way.
\end {proof}
We conclude this section by giving some simple result in the dilation theory of W*-dynamical systems, in particular we have the following proposition.
\begin{proposition}
Let $(\mathfrak{A},\Phi,\varphi)$ be a C*-dynamical system and $\Psi:\mathfrak{A}\rightarrow\mathfrak{A}$ an ucp-map with $\Phi(\Psi(a))=a$ for all $a\in\mathfrak{A}$. We have the following statement:
\begin{enumerate}
\item[1 - ] $\Psi(a)\in\mathcal{D}_{\Phi}$ for all $a\in\mathfrak{A}$, where $\mathcal{D}_{\Phi}$ is the multiplicative domain of the ucp-map $\Phi$;
\item[2 - ] $\Psi$ is a $\varphi$-adjoint of the ucp-map $\Phi$;
\item[3 - ] If $\varphi$ is faithful state, $\Psi$ is a homomorphism.
\end{enumerate}
\end{proposition}
\begin {proof}
1) For Kadison's inequality, for each $a\in\mathfrak{A}$ we have:
\[
0\leq\Phi(\Psi(a)^*\Psi(a))-\Phi(\Psi(a)^*)\Phi(\Psi(a))\leq \Phi(\Psi(a^*a))-a^*a=0,
\]
since $ \Psi(a^*)\Psi(a)\leq\Psi(a^*a)$. Then we can write that $\Psi(a)\in\mathcal{D}_{\Phi}$.
\newline
2) For each $a,b\in\mathfrak{A}$ we have:
\[
\varphi(a\Psi(b))=\varphi(\Phi(a\Psi(b))=\varphi(\Phi(a)\Phi(\Psi(b))=\varphi(\Phi(a)b),
\]
since $\Psi(b)\in\mathcal{D}_{\Phi}$. It follows that $\Phi^{\sharp}(b))=\Psi(b)$.
\newline
3) We recall that the multiplicative domain of a cp-map is a *-algebra, therefore for any $a\in\mathfrak{A}$ we obtain
\[
x=\Psi(a^*a)-\Psi(a)^*\Psi(a)\in\mathcal{D}_{\Phi}
\]
For the Kadison's inequality $x\geq 0$ and $\Psi(a)\in\mathcal{D}_{\Phi}$, it follows that
\[
\Phi(x)=\Phi(\Psi(a^*a))-\Phi(\Psi(a)^*\Psi(a))=\Phi(\Psi(a^*a))-\Phi(\Psi(a)^*)\Phi(\Psi(a)))=0.
\]
Then $x=0$ since $\varphi(\Phi(x))=\varphi(x)=0$ and $\varphi$ is a faithful state. 
\newline
We can write that $a\in\mathcal{D}_{\Psi}$ for all $a\in\mathfrak{A}$ where $\mathcal{D}_{\Psi}$ is the multiplicative domain of ucp-map $\Psi$. Then it is a homomorphism.
\end {proof}
We conclude this section with the following proposition:
\begin{corollary}
Let $(\mathfrak{M},\Phi,\varphi)$ be a W*-dynamical system. If there is an ucp-map $\Psi:\mathfrak{A}\rightarrow\mathfrak{A}$ such that $\Phi(\Psi(a))=a$ for all $a\in\mathfrak{A}$, the W*-dynamical system $(\mathfrak{M},\Phi,\varphi)$ admit a minimal reversible dilation.
\end{corollary}
\begin {proof}
From previous proposition the dynamics $\Phi$ admit as $\varphi$-adjoint the homomorphism $\Psi$, it follows that $(\mathfrak{M},\Psi,\varphi)$ admit a minimal reversible dilation. Therefore also our dynamical system admit a minimal reversible dilation.
\end {proof}

\section{Conclusion}
In this paper we have show that any dynamical system admit a covariant GNS as formuled by Niculescu, Str\"{o}h and Zsid\'{o} in \cite{NSZ}. Furthermore we have used this representation for to determine a reversible dilation for W*-dynamical system with multiplicative dynamics. In contrast to the existence of the CGNS representation, the W*-dynamical systems does not always have a reversible dilation (see \cite{Haag-Musat}) and in general it is not yet obvious when this can happen. A fundamental result in this direction is due to Haagerup and Musat in \cite{Haag-Musat}. They have proven that a W*-dynamical system $(\mathfrak{M},\Phi,\varphi)$ admits a reversible dilation if and only if the dynamics $\Phi:\mathfrak{M}\rightarrow\mathfrak{M}$ is a \emph{factorizable map} in the sense of Anantharaman-Delaroche in \cite{Antha}.
\newline
In finally, we observe that differently from dynamics $\Phi$, we have proved in Lemma \ref{lemma-fattor} that a linear contraction $\mathbf{U}_{\Phi,\varphi}$ associated to the dynamical system $(\mathfrak{M},\Phi,\varphi)$ is always factorizable through isometric operator on a Hilbert space (see \cite{Pisier}). 
 \section{Acknowledgments}
Thanks are due to Professor Laszlo Zsido (University of Rome Tor Vergata), for various discussions.

\end{document}